\def\EL{\mathop{\rm EL}}

\def\EEL{\mathop{\rm EEL}}

\def\area{\mathop{\rm area}}

\documentclass{article}

\usepackage{amsmath,amsthm}
\usepackage{url}
\usepackage{graphicx}
\usepackage{latexsym}
\usepackage{amsfonts}
\usepackage{amssymb}
\usepackage[font=sf]{caption}
\usepackage{subfigure}
\usepackage{verbatim}

\newtheorem{theorem}{Theorem}[section]
\newtheorem{lemma}[theorem]{Lemma}

\numberwithin{equation}{section}

\begin{document}

\title{Bounded Outdegree and Extremal Length on Discrete Riemann Surfaces}
\author{William E. Wood\footnote{Hendrix College, Conway, AR, billarama@gmail.com}}


\maketitle


\begin{abstract}
Let $T$ be a triangulation of a Riemann surface.  We show that the 1-skeleton of $T$ may be oriented so that there is a global bound on the outdegree of the vertices.  Our application  is to construct extremal metrics on triangulations formed from $T$ by attaching new edges and vertices  and subdividing its faces.  Such refinements provide a mechanism of convergence of the discrete triangulation to the classical surface.  We will prove a bound on the distortion of the discrete extremal lengths of path families on $T$ under the refinement process.  Our bound will depend only on the refinement and not on $T$. In particular, the result does not require bounded degree.
\end{abstract}

\section{Introduction}

Discrete conformal geometry is generally concerned with adapting ideas from classical conformal geometry to purely combinatorial objects.  A natural combinatorial analog of a Riemann surface is a triangulation.  In this paper, we discuss one of conformal geometry's most powerful tools, extremal length, in this combinatorial setting and establish bounds on how extremal length can change under the operation of combinatorial refinement.  Refinement is an important tool in connecting the discrete and classical settings as one expects classical results to appear as a limiting case of the discrete results under refinement.  Put another way, refinement is a mechanism through which discrete Riemann surfaces may better approximate their classical counterparts (c.f. \cite{rs}), establishing one possible model for a discrete surface and Teichm\"{u}ller theory.

Our main result is that under sufficient (and quite weak) regularity conditions, there is a global finite bound \emph{depending only on the refinement process} on the amount of distortion to the extremal length of a path family, extending results proved in the simply connected case to general surfaces \cite{wood}. This bound is entirely independent of the structure of the graph itself, meaning in particular that the degree of the graph is not an issue and thus holds for graphs with unbounded degree.

The reason degree restrictions may be avoided follows from the fact that there is a weaker notion of bounded degree satisfied by any graph on a Riemann surface, namely that it is always possible to orient the edges of a graph so that the outdegree is globally bounded.  This result is discussed in Section~\ref{sec:outdegree} and applications to extremal length are covered in Section~\ref{sec:extremallength}.  Some additional observations and extensions are remarked in our concluding Section~\ref{sec:conclusion}.

\section{Bounded Outdegree}\label{sec:outdegree}

Let $S$ be a Riemann surface with a locally finite simplicial complex $T$ consisting of vertices, edges, and triangular faces.    The graph formed by considering only the vertices and edges of $T$ is the \emph{complex graph} and is denoted $T^*$ (we disallow multigraphs and self-loops).  An \emph{orientation} of an edge is an ordering of the two vertices it bounds, commonly envisioned as an arrow pointing from a tail vertex to a head vertex.  An orientation on $T^*$ is an assignment of orientations to each of its edges, making $T^*$ into an oriented graph (not to be confused with the topological orientation on a manifold).  The \emph{outdegree} $\delta^+(v)$ of a vertex $v$ in an oriented graph $T^*$ is the number of edges for which $v$ is the tail, as opposed to \emph{degree} $\delta(v)$ which simply counts the number of edges that $v$ bounds (so $\delta^+(v)\leq \delta(v)$).

We are interested in finding orientations on $T^*$ for which $\delta^+$ is globally bounded, and we will show that this is always possible on a discrete Riemann surface with no reference to the graph itself.  Even if the graph has unbounded degree, we will prove the existence of an orientation such that $\delta^+(v)\leq 5$ for all $v$.  The key ingredient is the following observation of Chrobak and Eppstein \cite{eppstein}:

\begin{lemma} \label{th:outdeg3} The edges of a finite planar graph $G=(V,E)$ may be
oriented so that the outdegree of every interior vertex is at most three and the outdegree of each boundary vertex is at most two.
\end{lemma}

\begin{theorem}\label{outdegree}
Let $G$ be the graph formed from the edges and vertices of a simplicial complex $T$ on a Riemann surface $S$.  Then the edges of $G$ may be oriented so that the outdegree of every vertex is at most five.

\end{theorem}
\begin{proof} Topological classification of surfaces along with the pants decomposition theorem for hyperbolic surfaces (see e.g. \cite{teich}) allow us to assert the following:  there exists a collection of disjoint non-homotopic simple closed curves $\mathcal{B}$ on $S$ such that cutting along each curve in $\mathcal{B}$ will divide $S$ into a collection of surfaces $\mathcal{X}$ such that each $X\in \mathcal{X}$ is topologically one of the following: a plane, a torus, a half-infinite cylinder, or a sphere with at most three open disks removed.  Note that the last case is the usual pants decomposition of a hyperbolic surface into thrice-punctured spheres, but we make an additional cut around any open punctures; hence the inclusion of the half-infinite cylinder case.

Since $G$ is a simplicial complex, each simple closed curve in $\mathcal{B}$ is homotopic to a cycle subgraph of $G$ as embedded in $S$.  Let $\mathcal{C}$ be the collection of such cycle graphs.  We assume that if two cycles in $\mathcal{C}$ intersect at distinct vertices, then both cycles share an edge path connecting those two vertices.  We now cut  along these cycle graphs to obtain a collection $\mathcal{G}$ of subgraphs of $G$, each of which defines a simplicial complex on one of the surfaces described above.  A cycle in $\mathcal{C}$ now corresponds to two cycles in the members of $\mathcal{G}$; otherwise, edges and vertices of the members of $\mathcal{G}$ are in one-to-one correspondence with those of $G$.

Our task now is to extend Lemma~\ref{th:outdeg3} to the surfaces in our list and then show how to preserve the bounded outdegree property when we reassemble $G$.

In the plane, we begin by choosing a face $T_0$ in the complex and let $G_0$ be the corresponding three-edge cycle graph.  Now let $T_1$ be any disk subcomplex in $T$ containing $T_0$ in its interior and let $A_1 = T_1 \setminus T_0$ be the complex for the annulus formed by removing the face $T_0$ from $T_1$.  Similarly, let $T_2$ be a subcomplex of a disk in $T$ containing $T_2$ in its interior and define a new annulus $A_2$ by removing $A_1$ and $A_0$ from $T_2$.  Continuing in this fashion, we construct a sequence $A_1, A_2, A_3,\ldots$ of annular complexes that fill the plane and intersect only on their boundaries (except $A_1$, which meets the face $T_1$; it is thus convenient to introduce $A_0=T_0$ as a degenerate annulus).

Now apply Lemma~\ref{th:outdeg3} to orient each of the $A_i^*$ so that each interior vertex has outdegree at most three and the boundary cycles have outdegree two.  We now reorient the boundary cycles so that progression from tail to head progresses around the cycle in a counterclockwise direction.  (Counterclockwise here is inherited from how the annuli $A_i$ sit in the plane with an arbitrarily chosen topological orientation.)  This new orientation on the $A_i^*$ may increase the outdegree of a boundary vertex by one, so the outdegree of any vertex is at most three.

We now orient $G$ in the obvious way, assigning an edge the orientation it inherits from the annuli in which it lies.  Note that this is well-defined for an edge lying on the intersection of two of the $A_i$ because the orientations on boundary cycles were chosen consistently.  This orientation on $G$ gives a worst-case outdegree of five on a boundary vertex: three from each of the annuli it bounds minus one for double-counting the outgoing edge on the boundary cycle. This proves the theorem for the plane.

The results for the remaining component surfaces follow immediately, as any complex graph of a half-open cylinder or a subset of the sphere is already planar.  For the torus, we simply cut along any edge cycle that turns $G$ into a graph complex $G'$ of a finite cylinder.  $G'$ is a finite planar graph and we assign  it the orientation guaranteed by Lemma~\ref{th:outdeg3}.  We then apply the same reorientation and regluing method we used for the plane to induce an orientation with outdegree bounded by five.
\end{proof}

\section{Discrete Extremal Length}\label{sec:extremallength}

Our application of the bounded outdegree theorem is to construct extremal metrics on discrete Riemann surfaces.

Let  $G=(V,E)$ be a graph and $\Gamma$ a non-empty collection of
finite or infinite vertex paths in $G$.  A \emph{metric} on $G$
is a function $m:V\to [0,\infty)$. The value $m(x)$ is the
\emph{$m$-weight} or \emph{$m$-measure} of $x$. The \emph{area} of
$m$ is
$$\mathrm{area}(m) = \sum_{v\in V} m(v)^2,$$ and a metric is called \emph{admissible} if its area is  finite and non-zero.  Let
${\mathcal{M}(G)}=\{m:0<\area(m) <\infty\}$ be the set of admissible
metrics on $X$.

A \emph{(vertex) path} in $G$ is a sequence of vertices in $V$ such that consecutive vertices are either adjacent or identical.  For a path $\gamma=\{a_0,a_1,\ldots\}\subset V$ and $m\in{\mathcal{M}(G)}$, define
the $m$-\emph{length} to be $L_m(\gamma) = \sum_{j=1}^\infty m(a_j) = \sum_{v\in \gamma}
m(v)$. For a collection $\Gamma$ of paths
in $G$, which we call a \emph{path family}, define $\displaystyle L_m(\Gamma) = \inf_{\gamma\in \Gamma} L_m(\gamma)$ and the
\emph{extremal length}
$$\EL(\Gamma) = \sup_{m\in \mathcal{M}(G)} \left\{
\frac{L_m(\Gamma)^2}{\area(m)}\right\}.$$  An \emph{extremal metric} $\mu$ for $\Gamma$ is one that realizes the extremal length, i.e. $\EL(\Gamma)=L_\mu(\Gamma)^2 / \area(\mu)$.  Such metrics always exist by \cite{cannonacta},\cite{heschrammhp}.  By rescaling, we may assume an extremal metric has unit area.


\begin{figure}
\begin{center}
\includegraphics[width=3in]{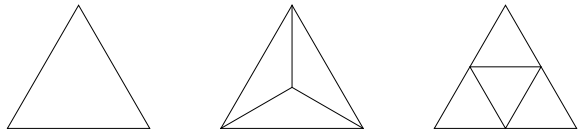}

\end{center}
\caption{Examples of refinement, left to right:  the identity, barycentric, and hexagonal refinements.  All are $b$-bounded with $b=0,0,1$, respectively.  They are also strongly bounded with $c=0,1,2$.}\label{refineexamples}
\end{figure}

We will relate the extremal length of a path family in $G$ to the extremal length of a corresponding path family in a subdivision of $G$.  A \emph{refinement} $rG=(rV, rE)$ of $G$ is a triangulation with an injection $\iota: V\to rV$ and a mapping of each edge in $(a,b)\in E$ to a finite path in $rG$ from $\iota(a)$ to $\iota(b)$ such that no other vertices in the path lie in the image of $\iota$ or in another path.  The vertices in $rV$ that form the path from $\iota(a)$ to $\iota(b)$, excluding the endpoints, are said to be \emph{attached} to the edge $(a,b)$.  We will often suppress mention of the mapping $\iota$ and consider $V\subset rV$. The refinement is \emph{$b$-bounded} if the maximum length of any such path is $b+2$ (the 2 is to discount the endpoints).  Some examples are shown in Figure~\ref{refineexamples}.

 We will require some consistency conditions on our path families.  We say a path family $\Gamma$ is \emph{regular} if it satisfies the following properties.  If for any triangle $T$ in $G$ with vertices $a,b,c$ and the property that each pair $a b$, $b c$, and $a c$ appear in some curves in $\Gamma$ (not necessarily all the same curve), then for any $\gamma\in\Gamma$  of the form $\gamma = x a b y$, where $x$ and $y$ are sequences of vertices, we have $x a c b y \in \Gamma$, and similarly for the remaining pairs $b c$ and $a c$.  We then say $\Gamma$ \emph{surrounds} $T$.  The other property we require for regularity is that if $\gamma_1 = x_1 a y_1$ and $\gamma_2 = x_2 a y_2$ are paths in $\Gamma$ where  $x_1, x_2, y_1, y_2$ are sequences of vertices and $a\in V$, then the path $x_1 a y_2$ is also in $\Gamma$.

 Let $\Gamma$ be a regular curve family on a triangulation graph $G$ and let $rG$ be a refinement of $G$.  We define a corresponding curve  family $r\Gamma$ on $rG$ as follows.  A path $\gamma\in rG$ is in $r\Gamma$ if $\gamma$ can be expressed as a sequence of finite vertex paths $\gamma=\alpha_0 \alpha_1 \ldots$ such that for each $\alpha_i$ there is a triangle $T_i$ surrounded by $\Gamma$ for which each vertex in the path $\alpha_i$ is either inside or attached to $T_i$.  The point of regularity is to ensure that if a subpath in $\Gamma$ reaches a vertex of a triangle that intersects some refined path then we may extend that subpath by either of the other vertices of the triangle and remain a subpath of a path in $\Gamma$.

The motivating example for these definitions is a path family defined as the set of paths connecting two simple closed curves or ends of the discrete surface.  The extremal length of an annulus, for example, is the extremal length of the family of curves connecting its boundary cycles.  The refinement of such a family will be the set of curves connecting the refined boundary cycles, i.e. the refined edges.  A limiting case in the simply connected case is the type problem, which asks whether the extremal length of a family of curves connecting a base point to infinity (a limiting annulus) is finite.

Our main result is the following.
\begin{theorem}\label{th:mainth} Let $G$ be a triangulation of a discrete Riemann surface, $r G$ a $b$-bounded refinement of $G$, and $\Gamma$ a regular curve family.  Then there exists a $k\geq 1$ depending only on $b$ such that $\frac{1}{k}\EL(r\Gamma) \leq \EL(\Gamma) \leq k\EL(r\Gamma)$.  In particular, $\EL(\Gamma)$ is finite if and only if $\EL(r\Gamma)$ is.
\end{theorem}

\begin{proof} Our approach is to use an extremal metric on one graph to construct one on the other.

Suppose $\mu$ is a unit-area extremal metric on $rG$, so $\displaystyle{\area(\mu) =\sum_{v\in V} \mu(v)^2 = 1}$ and $\displaystyle{\EL(rG) = \inf_{\gamma\in\Gamma} L_\mu(\gamma)^2}$.

Let $\theta^\star(v)=\displaystyle{\max_{w\in\epsilon(v)} \mu(w)}$, where $\epsilon(v)$ is the set of  all vertices $w\in rV$ that are attached to an edge bounded by $v$.  Define a metric $\theta$ on $G$ by $\theta(v) = \max \{\theta^\star(v), \mu(v)\}$ (recall our notational abuse that $\mu(v)=\mu(\iota(v))$).

Then $$\area(\theta) = \sum_{v\in V} \theta(v)^2  \leq 2 \sum_{w\in rV} \theta(w)^2 =2 \area(\theta) = 2.$$
The inequality above comes from the fact that a vertex in $rV$ is contained in at most two of the $\epsilon(v)$ because an edge is bounded by exactly two vertices.  In other words, no vertex in $rV$ can contribute its $\mu$-weight to $\theta$ more than twice.

Now let $\gamma\in\Gamma$ and define $\gamma_r \in r\Gamma$ to be the vertex path formed by traveling along $\gamma$ as it sits in $rG$.  More precisely, let $E(\gamma)$ denote the collection of edges $v w\in G$ such that $v w$ appears in the vertex sequence defining $\gamma$.  To define $\gamma_r$, if  $v w\in E(\gamma)$ and $\iota(v) x_1 x_2 \ldots x_j \iota(w)$ is the refined edge in $rG$, we replace $v w$  with $\iota(v) x_1 x_2 \ldots x_j \iota(w)$ to form $\gamma_r$.  It is also convenient to define $r(e)$ to be the set of refined vertices $x_1\ldots x_j$.  We now compare the lengths of these paths in their respective metrics.

$$L_\mu(\gamma_r) = \sum_{z\in \gamma_r} \mu(z) = \sum_{v\in\gamma} \mu(v)+\sum_{e\in E(\gamma)} \sum_{z\in r(e)} \mu(z)$$
$$ \leq L_\theta(\gamma) + \sum_{e\in E(\gamma)} b \max_{z\in r(e)} \theta(z) \leq L_\theta(\gamma)+ \sum_{v\in \gamma} b\theta(v)  = (1+b)L_\theta(\gamma)$$

We have thus shown that every path $\gamma\in\Gamma$ has a ``shadow path" $\gamma_r \in r\Gamma$ whose length, scaled by the fixed constant $b+1$, is shorter that $\gamma$.  Carrying this correspondence to the infimum, it follows that  $L_\mu(r\Gamma)\leq (b+1)L_\theta (\Gamma)$ and thus $$\EL(r\Gamma) = L_\mu(r\Gamma)^2 \leq (b+1)^2 L_\theta(\Gamma)^2 = 2\cdot \frac{ (b+1)^2 L_\theta(r\Gamma)^2}{2}$$
$$\leq 2 \frac{ (b+1)^2 L_\theta(\Gamma)^2}{\area(\theta)} \leq 2(b+1)^2\cdot\sup_{m\in \mathcal{M}(G)} \frac{ L_m(\Gamma)^2}{\area(m)}=(b+1)^2\EL(\Gamma)$$
proving half of our theorem.

The reverse direction will be approached similarly, but our metric construction will require our results from Section~\ref{sec:outdegree} to work without degree requirements on $G$.

Let $\theta$ be an extremal metric on $G$.  Apply Theorem~\ref{outdegree} to $G$ to orient the edges so that $\delta^+(v)\leq 5$ for every $v\in V$.  Define a metric $\mu$ on $rG$ by
\begin{displaymath}\mu(v) = \left\{
  \begin{array}{ll}
     \theta(v) & \textrm{if\ } v\in V
     \\ \\
     \theta(w) & \mathrm{if\ } v\in r(e) \mathrm{\ and\ }  w \\
     & \mathrm{\ is\ the\ tail\ of\ the\ oriented\ edge\ } e
   \\ \\
     0 & \textrm{otherwise.}
  \end{array}
  \right.
\end{displaymath}
The finite area condition follows easily from the bounded outdegree condition because no vertex can be the tail of more than five vertices:
$$
\area(\mu) = \sum_{v\in rV} \mu(v)^2 = \sum_{v\in V} \mu(v)^2 +
\sum_{v\in rV\setminus V} \mu(v)^2 $$ $$\leq \sum_{v\in V} \mu(v)^2
+ \sum_{v\in V} 5b \mu(v)^2
 =(1+5b) \area(\mu)=1+5b.$$

Let $\gamma_r \in r\Gamma$.  We will construct a path $\gamma\in\Gamma$ such that $L_\mu(\gamma_r)= L_\theta(\gamma)$, and the desired result will carry over to extremal length just as in the first part of the proof.

\begin{figure}
\begin{center}
\includegraphics[width=1.5in]{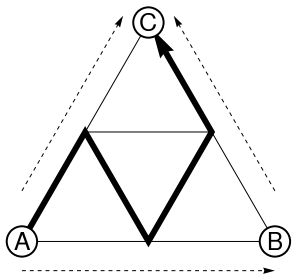}\end{center}
\caption{Shadow path construction.  $ABC$ is a triangle in $G$, and the bold edges indicate a five-vertex path through the hexagonal refinement.  If $G$ is oriented according to the dashed arrows, then the corresponding shadow path is $AAABC$.}\label{shadowpath}
\end{figure}

Construct $\gamma$ inductively.  Assume for convenience that the starting vertex of $\gamma_r$ is $p\in V$, so that our first iteration is $\gamma^0 =\gamma_r^0 = p$ (after we illustrate the construction, the reader may wish to verify that this assumption is not a problem). Assume $\gamma^s$ has been constructed so that $L_\mu(\gamma_r^s)=L_\theta(\gamma^s)$ for a subpath $\gamma^s$ of $\gamma$.  Let $v$ be the final point of $\gamma^s$ and $v'$ the final point of $\gamma_r^s$.  Assume that if $v'\in V$ then $v'=v$ and if $v'\in r(e)$ for some $e\in E$ then $v$ is an endpoint of $e$.  Our approach is to extend $\gamma^s$ by a single vertex $w$ to form $\gamma^{s+1}$ so that these properties are preserved and so that $L_\mu(\gamma_r^{s+1}) = L_\theta(\gamma^{s+1})$ for a subpath $\gamma_r^{s+1}$ of $\gamma_r$ .

Let $w'$ be the vertex following $v'$ in $\gamma_r$.  There are three possibilities for how to  extend $\gamma^s$ to $\gamma^{s+1}$ by a vertex $w$.  If $w'$ is not in $r(e)$ for any edge $e$, then set $\gamma^{s+1} = \gamma^s$.  We thus maintain $L_\mu(\gamma_r^{s+1}) = L_\theta(\gamma^{s+1})$ because $L_\mu (w')=0$.

If $w' \in V$, then set $w=w'$.  The fact that $G$ is a triangulation guarantees that $v$ is equal or adjacent to $w$ and thus $\gamma^{s+1}$ is still a path, and this path is in $\Gamma$ because $r$ is a refinement and $\Gamma$ is regular.

Otherwise, we have that $w'\notin V$ but $w'\in r(e)$ for some edge $e$.  We adopt the rule \emph{always move to the base of the arrow}, meaning we take $w$ to be the tail of the directed edge $e$.  We again have $L_\mu(\gamma_r^{s+1}) = L_\theta(\gamma^{s+1})$ by the definition of $\mu$ and the same guarantee that $\gamma^{s+1}$ is a path in $\Gamma$.  See Figure~\ref{shadowpath}.

This construction guarantees $L_\theta(\Gamma) \leq L_\mu(r\Gamma)$ because each path in $r\Gamma$ has a shadow of equal length in $\Gamma$, so the infimum of all path lengths in $r\Gamma$ cannot be smaller than in $\Gamma$. We thus obtain
$$\EL(\Gamma) = L_\theta(\Gamma)^2 \leq  L_\mu(r\Gamma)^2 =(1+5b)\cdot \frac{ L_\mu(r\Gamma)^2}{1+5b}$$
$$\leq (1+5b) \frac{ L_\mu(r\Gamma)^2}{\area(\mu)} \leq (1+5b)\cdot\sup_{m\in \mathcal{M}(rG)} \frac{ L_m(\Gamma)^2}{\area(m)}=(1+5b)\EL(r\Gamma)$$
\end{proof}

\section{Extensions and Remarks}\label{sec:conclusion}

The methods of this paper may be generalized to several related questions about discrete extremal length.  The results we mention in this section are proved for the simply connected case in \cite{wood} and the reader is invited to use the tools we have developed to extend those methods to discrete Riemann surfaces.

There is another way to define discrete extremal length.  Instead of defining metrics on the vertices of $G$, we could just as easily have defined it on the edges.  This gives an alternate definition of extremal length called \emph{edge extremal length} (EEL), with paths, lengths, and regularity all defined analogously to the vertex case.   In the plane, edge extremal length captures the behavior of random walks and electric networks whereas the vertex extremal length we have studied relates to circle packings. They are not equivalent for graphs of unbounded degree. See  \cite{doyle},\cite{duffin},\cite{heschrammhp}.  We get something similar to Theorem~\ref{th:mainth} for edge extremal length, but with a stronger refinement condition.

A refinement is \emph{$(b,c)$-strongly bounded} if it is $b$-bounded and for any cell in $F$, any edge $e\in E$, and any vertex $w\in rV$ bounding $e$ or in $r(e)$, we have that $w$ is the endpoint of at most $c$ edges in $rE$ whose other endpoint lies in $F$. The proof requires a similar strategy to that used in the vertex case, but bounded outdegree is not required.

\begin{theorem} Let $G$ be a triangulation of a discrete Riemann surface, $r G$ a $(b,c)$-strongly bounded refinement of $G$, and $\Gamma$ a regular family of edge paths.  Then there exists a $k\geq 1$ depending only on $b$ and $c$ such that $\frac{1}{k}\EEL(r\Gamma) \leq \EEL(\Gamma) \leq k\EEL(r\Gamma)$.  In particular, $\EEL(\Gamma)$ is finite if and only if $\EEL(r\Gamma)$ is.
\end{theorem}

It is also worth noting that similar results can be obtained when $G$ is a polygonal complex and not necessarily a triangulation.  Recall in the proof of Theorem~\ref{th:mainth} that the metric $\mu$ was constructed by keeping track of a path as it moved through the triangles.  If the cells are not triangles, then vertices on the cells may not be adjacent.  A similar construction works, but the maximal number of edges on the polygonal cells appears in the bound.

\bibliographystyle{amsplain}
\bibliography{extremalriemann}

\providecommand{\bysame}{\leavevmode\hbox to3em{\hrulefill}\thinspace}
\providecommand{\MR}{\relax\ifhmode\unskip\space\fi MR }
\providecommand{\MRhref}[2]{%
  \href{http://www.ams.org/mathscinet-getitem?mr=#1}{#2}
}
\providecommand{\href}[2]{#2}
\begin{thebibliography}{1}

\bibitem{cannonacta}
J.~W. Cannon, \emph{The combinatorial {Riemann} mapping theorem}, Acta
  Mathematica \textbf{173} (1994), 155--234.

\bibitem{eppstein}
M.~Chrobak and D.~Eppstein, \emph{{Planar orientations with low out-degree and
  compaction of adjacency matrices}}, Theoretical Computer Science \textbf{86}
  (1991), no.~2, 243--266.

\bibitem{doyle}
P.~G. Doyle and J.~L. Snell, \emph{Random walks and electric networks}, The
  Carus Mathematical Monographs, no.~22, Math. Association of America, 1984.

\bibitem{duffin}
R.~J. Duffin, \emph{The extremal length of a network}, Journal of Mathematical
  Analysis and Applications \textbf{5} (1962), 200--215.

\bibitem{heschrammhp}
Z.-X. He and O.~Schramm, \emph{Hyperbolic and parabolic packings}, Discrete \&
  Computational Geom. \textbf{14} (1995), 123--149.

\bibitem{teich}
J.~H. Hubbard, \emph{Teichm\"{u}ller theory and application to geometry,
  topology, and dynamics, volume 1}, Matrix Editions, 2006.

\bibitem{rs}
B.~Rodin and D.~Sullivan, \emph{The convergence of circle packings to the
  {Riemann} mapping}, J. Differential Geometry \textbf{26} (1987), 349--360.

\bibitem{wood}
W.~E. Wood, \emph{Combinatorial modulus and type of refined graphs}, Topology
  and its Applications \textbf{156} (2009), no.~17, 2747--2761,
  doi:10.1016/j.topol.2009.02.013.

\end{thebibliography}

\end{document}